# BASES OF SOLUTIONS FOR LINEAR CONGRUENCES

Florentin Smarandache
University of New Mexico
200 College Road
Gallup, NM 87301, USA
E-mail: smarand@unm.edu

In this article we establish some properties regarding the solutions of a linear congruence, bases of solutions of a linear congruence, and the finding of other solutions starting from these bases.
This article is a continuation of my article "On linear congruences".

### §1. Introductory Notions

**Definition 1.** (linear congruence)
We call linear congruence with $n$ unknowns a congruence of the following form:
$$a_1 x_1 + ... + a_n x_n \equiv b (\bmod m) \qquad (1)$$
where $a_1,...,a_n, m \in \mathbb{Z}, n \geq 1$, and $x_i, i = \overline{1,n}$, are the unknowns.
The following theorems are known:
**Theorem 1.** The linear congruence (1) has solutions if and only if $(a_1,...,a_n,m,b) | b$.

**Theorem 2.** If the linear congruence (1) has solutions, then: $|d| \cdot |m|^{n-1}$ is its number of distinct solutions. (See the article "On the linear congruences".)

**Definition 2.** Two solutions $X = (x_1,...,x_n)$ and $Y = (y_1,...,y_n)$ of the linear congruence (1) are distinct (different) if $\exists i \in \overline{1,n}$ such that $x_i \not\equiv y_i (\bmod m)$.

### §2. Definitions and proprieties of congruences

We'll present some arithmetic properties, which will be used later.
**Lemma 1.** If $a_1,...,a_n \in \mathbb{Z}$, $m \in \mathbb{Z}$, then:
$$\frac{(a_1,...,a_n,m) \cdot m^{n-1}}{(a_1,m) \cdot ... \cdot (a_n,m)} \in \mathbb{Z}$$
The proof is done using complete induction for $n \in \mathbb{N}^*$.
When $n = 1$ it is evident.
Considering that it is true for values smaller or equal to $n$, let's proof that it is true for $n+1$.
Let's note $x = (a_1,...,a_n)$. Then:

$(a_1,...,a_n,a_{n+1},m) \cdot m^n = \left[(x,a_{n+1},m) \cdot m^{2-1}\right] \cdot m^{n-1}$, which, in accordance to the induction hypothesis, is divisible by:

$\left[(x,m) \cdot (a_{n+1},m)\right] \cdot m^{n-1} = \left[(a_1,...,a_n,m) \cdot (a_{n+1},m)\right] \cdot m^{n-1} = \left[(a_1,...,a_n,m) \cdot m^{n-1}\right] \cdot (a_{n+1},m)$,

which is divisible, also in accordance with the induction hypothesis, by

$\left[(a_1,m) \cdot ... \cdot (a_n,m)\right] \cdot (a_{n+1},m) = (a_1,m) \cdot ... \cdot (a_n,m) \cdot (a_{n+1},m)$.

**Theorem 3**. If $X^0$ constitutes a (particular) solution of the linear congruence (1), and $p = \prod_{i=1}^{n}(a_i,m)$, then:

$$X_i \equiv x_i^0 + \frac{m}{(a_i,m)} t_i, \quad 0 \leq t_i < (a_i,m), \quad t_i \in \mathbb{N} \quad (*)$$

($i$ taking values from 1 to $n$) constitute $p$ distinct solutions of (1).

*Proof:*

Because the module of the congruence (m) is sub-understood, we omitted it, and we will continue to omit it.

$$\sum_{i=1}^{n} a_i x_i = \sum_{i=1}^{n} a_i x_i^0 + \sum_{i=1}^{n} \frac{a_i m}{(a_i,m)} t_i \equiv b + 0,$$ therefore there are solutions. Let's show that they are also distinct.

$$x_i^0 + \frac{m}{(a_i,m)} \alpha \not\equiv x_i^0 + \frac{m}{(a_i,m)} \beta, \quad \text{for} \quad \alpha, \beta \in \mathbb{N}, \quad \alpha \neq \beta, \text{ and } 0 \leq \alpha, \beta < (a_i,m),$$

because the set:

$$\left\{\frac{m}{(a_i,m)} t_i \mid 0 \leq t_i < (a_i,m), \quad t_i \in \mathbb{N}\right\} \subseteq \{0,1,...,n-1\},$$ which constitutes a complete system of residues modulo $m$, and $\frac{m}{(a_i,m)}\alpha \neq \frac{m}{(a_i,m)}\beta$, for $\alpha$ and $\beta$ previously defined.

Therefore the theorem is proved.

<p align="center">*<br>* *</p>

One considers the $Z$-module $A$ generated by the vectors $V_i$, where

$$V_i^* = \left(\underbrace{0,...,0}_{i-1 \text{ times}}, \frac{m}{(a_i,m)}, \underbrace{0,...,0}_{n-i \text{ times}}\right), i = \overline{1,n}, \text{ from } \mathbb{Z}^n. \text{ The module } A \text{ has the rank } n, (n \geq 1).$$

We could note it $A = \{v_1,...,v_n\}$.

We'll introduce a few new terms.

**Definition 3.** Two solutions (vectors solution) $X$ and $Y$ of congruence (1) are called independent if $X - Y \notin A$. Otherwise, they are called dependent solutions.

**Remark 1.** In other words, if $X$ is a solution of the congruence (1), then the solution $Y$ of the same congruence is independent of $X$, if it was not obtained from $X$ by applying the formula (*) for certain values of the parameters $t_1, ..., t_n$.

**Definition 4.** The solutions $X^1, ..., X^n$ are called **independent (all together)** if they are independent two by two.
Otherwise, they are called **dependent solutions** (**all together**).

**Definition 5.** The solutions $X^1, ..., X^n$ of the congruence (1) constitute a base for this congruence, if $X^1, ..., X^n$ are independent amongst them, and with their help one obtains all (distinct) solutions of the congruence with the procedure (*) using the parameters $t_1, ..., t_n$.

**Some proprieties of the linear congruences solutions:**
1) If the solution $X^1$ is independent with the solution $X^2$ then $X^2$ is independent with $X^1$ (the commutative property of the relation "independent").
2) $X^1$ is not independent with $X^1$.
3) If $X^1$ is independent with $X^2$, $X^2$ is independent with $X^3$, it does not imply that $X^1$ is independent with $X^3$ (the relation is not transitive).
4) If $X$ is independent with $Y$, then $X$ is independent with $Y$.
   Indeed, if $Y$ is dependent with $Y$, then $X - Y = \underbrace{(X - Y)}_{\notin A} + \underbrace{(Y - Y_1)}_{\in A} = Z$.

If $Z \in A$, it results that $(X - Y) = Z - (Y - Y_1) \in A$ because $A$ is a $Z$-module. Absurdity.

$$*$$
$$*\quad *$$

**Theorem 4.** Let's note $P_1 = (a_1, ..., a_n, m) \cdot |m|^{n-1}$ and $P_2 = (a_1, m) \cdot ... \cdot (a_n, m)$ then the linear congruence (1) has the base formed of: $\dfrac{P_1}{P_2}$ solutions.

*Proof:*

$P_1 > 0$ and $P_2 > 0$, from Lemma 1 we have $\dfrac{P_1}{P_2} \in \mathbb{N}^*$, therefore the theorem has sense (we consider LCD as a positive number).
$P_1$ represents the number of distinct solutions (in total) of congruence (1), in accordance to theorem 2.
$P_2$ represents the number of distinct solutions obtained for congruence (1) by applying the procedure (*) (allocating to parameters $t_1, ..., t_n$ all possible values) to a single particular solution.

Therefore we must apply the procedure (*) $\dfrac{P_1}{P_2}$ times to obtain all solutions of the congruence, that is, it is necessary of exact $\dfrac{P_1}{P_2}$ independent particular solutions of the congruence. That is, the base has $\dfrac{P_1}{P_2}$ solutions.

**Remark 2.** Any base of solutions (for the same linear congruence) has the same number of vectors.

### §3. Method of solving the linear congruences

In this paragraph we will utilize the results obtained in the precedent paragraphs.
Let's consider the linear congruence (1) with $(a_1,...,a_n,m) = d \mid b$, $m \neq 0$.

- we determine the number of distinct solutions of the congruence: $P_1 = |d| \cdot |m|^{n-1}$;

- we determine the number of solutions from the base: $S = \dfrac{P_1}{\prod_{i=1}^{n}(a_i,m)}$;

- we construct the $Z$-module $A = \{V_1,...,V_n\}$, where
$$V_i^t = \left(\underbrace{0,...,0}_{i-1\ times},\ \dfrac{m}{(a_i,m)},\ \underbrace{0,...,0}_{n-i\ times}\right),\ i = \overline{1,n};$$

- we search to find $s$ independent (particular) solutions of the congruence;
- we apply the procedure (*) as follows:

if $X^j$, $j = \overline{1,s}$, are the $s$ independent solutions from the base, it results that
$$X^{j(t_1,...,t_n)} = \left(x_i^j + \dfrac{m}{(a_i,m)} t_i\right),\ i = \overline{1,n}, \qquad (*)$$
are all $P_1$ solutions of the linear congruence (1),
$$j = \overline{1,s},\ t_1 \times ... \times t_n \in \{0,1,2,...,d_1-1\} \times ... \times \{0,1,2,...,d_n-1\},$$
where $d_i = |(a_i,m)|$, $i = \overline{1,n}$.

**Remark 3.** The correctness of this method results from the anterior paragraphs.

**Application.** Let's consider the linear non-homogeneous congruence $2x - 6y \equiv 2 \pmod{12}$. It has $(2,6,12) \cdot 12^{2-1} = 24$ distinct solutions. Its base will have $24 : [(2,12) \cdot (6,12)] = 2$ solutions.

$$V_1^t = (6,0),\ V_2^t = (0,2)\ \text{and}\ A = \{V_1,V_2\} = \{(6t_1, 2t_2)^t \mid t_1, t_2 \in \mathbb{Z}\}.$$

The solutions $x \equiv 7 \pmod{12}$ and $y \equiv 4 \pmod{12}$, $x \equiv 1$ and $y \equiv 0$ are dependent because:
$$\binom{7}{0} - \binom{1}{0} = \binom{6}{4} = 1\binom{6}{0} + 2\binom{0}{2} \in A.$$

But $\binom{4}{1}$ is independent with $\binom{0}{1}$ because $\binom{4}{1} - \binom{0}{1} \notin A$.

Therefore, the 24 solutions of the congruence can be obtained from:
$$\begin{cases} x \equiv 1 + 6t_1, & 0 \leq t_1 < 2, \ t_1 \in \mathbb{N} \\ y \equiv 0 + 2t_2, & 0 \leq t_2 < 6, \ t_2 \in \mathbb{N} \end{cases}$$
and
$$\begin{cases} x \equiv 4 + 6t_1, & 0 \leq t_1 < 2, \ t_1 \in \mathbb{N} \\ y \equiv 1 + 2t_2, & 0 \leq t_2 < 6, \ t_2 \in \mathbb{N} \end{cases}$$
by the parameterization $(t_1, t_2) \in \{0,1\} \times \{0,1,2,3,4,5\}$.

$$\begin{cases} x \equiv 1 + 6t_1 \\ y \equiv 0 + 2t_2 \end{cases} \Rightarrow \binom{1}{0}, \binom{1}{2}, \binom{1}{4}, \binom{1}{6}, \binom{1}{8}, \binom{1}{10}, \binom{7}{0}, \binom{7}{2}, \binom{7}{4}, \binom{7}{6}, \binom{7}{8}, \binom{7}{10}.$$

$$\begin{cases} x \equiv 4 + 6t_1 \\ y \equiv 1 + 2t_2 \end{cases} \Rightarrow \binom{4}{1}, \binom{4}{3}, \binom{4}{5}, \binom{4}{7}, \binom{4}{9}, \binom{4}{11}, \binom{10}{1}, \binom{10}{3}, \binom{10}{5}, \binom{10}{7}, \binom{10}{9}, \binom{10}{11};$$

which constitute all 24 distinct solutions of the given congruence; $\binom{0}{1}$ means: $x \equiv 1 \pmod{12}$ and $y \equiv 0 \pmod{12}$; etc.

### REFERENCES

[1] C. P. Popovici – "Teoria numerelor", Editura didactică și pedagogică, București, 1973.

[2] F. Gh. Smarandache – "Rezolvarea ecuațiilor și a sistemelor de ecuații liniare în numere întregi", Lucrare de licență, Universitatea din Craiova, 1979.

[Published in "Bulet. Univ. Brașov", series C, vol. XXII, pp. 25-31, 1980; and in "Bulet. Șt. și Tehn. al Instit. Polit. "Traian Vuia", Timișoara", fasciculul 2, tomul 26 (40), pp. 13-6; MR: 83c: 10024]